\documentclass[a4paper,oneside]{amsart}
  
\usepackage{amssymb,pstricks,amscd}
\usepackage{setspace}
\usepackage[totalwidth=14cm,totalheight=20cm]{geometry}
\usepackage{graphicx}

\newcommand{\cT}{{\mathcal T}}

\newcommand{\cGH}{{\mathcal G\mathcal H}}

\newcommand{\C}{{\mathbb C}}
\newcommand{\N}{{\mathbb N}}
\newcommand{\R}{{\mathbb R}}

\newcommand{\SL}{\mathrm{SL}}
\newcommand{\PSL}{\mathrm{PSL}}
\newcommand{\SO}{\mathrm{SO}}
\newcommand{\OO}{\mathrm{O}}
\newcommand{\U}{\mathrm{U}}
\newcommand{\Hom}{\mathrm{Hom}}
\newcommand{\AdS}{\mathrm{AdS}}
\newcommand{\T}{\mathcal{T}}
\newcommand{\HH}{\mathbb{H}}

\newtheorem{prop}{Proposition}[section]

\newtheorem{remark}[prop]{Remark}

\newtheorem{question}[prop]{Question}

\begin{document}

\title{Some open questions on anti-de Sitter geometry}

\author[]{Thierry Barbot}
\address{TB: Laboratoire d'analyse non lin\'eaire et g\'eom\'etrie,
Universit\'e d'Avignon et des pays de Vaucluse,
33 rue Louis Pasteur, 84 018 Avignon, France}
\email{thierry.barbot@univ-avignon.fr}
\author[]{Francesco Bonsante}
\address{FB: Universit\`a degli Studi di Pavia,
Via Ferrata, 1, 27100 Pavia, Italy}
\email{francesco.bonsante@unipv.it}
\author[]{Jeffrey Danciger}
\address{JD: University of Texas-Austin, Department of Mathematics, 
1 University Station C1200, Austin, TX 78712-0257, USA}
\email{jdanciger@math.utexas.edu}
\author{William M. Goldman}
\address{WMG: Department of Mathematics, University of Maryland,
College Park, MD 20742, USA }
\email{wmg@math.umd.edu}
\author[]{Fran\c{c}ois Gu\'eritaud}
\address{FG: CNRS and Laboratoire Paul Painlev\'e, Universit\'e Lille 1,
59 655 Villeneuve d'Ascq Cedex, France}
\email{francois.gueritaud@math.univ-lille1.fr}
\author[]{Fanny Kassel}
\address{FK: CNRS and Laboratoire Paul Painlev\'e, 
Universit\'e Lille 1, 
59 655 Villeneuve d'Ascq Cedex, France}
\email{fanny.kassel@math.univ-lille1.fr}
\author[]{Kirill Krasnov}
\address{KK: School of Mathematical Sciences,
University of Nottingham, Nottingham, NG7 2RD, UK}
\email{kirill.krasnov@nottingham.ac.uk}
\author[]{Jean-Marc Schlenker}
\address{JMS: Institut de Math\'ematiques de Toulouse, UMR CNRS 5219,
Universit\'e Paul Sabatier, 31062 Toulouse Cedex~9, France}
\email{schlenker@math.univ-toulouse.fr}
\author[]{Abdelghani Zeghib}
\address{AZ: CNRS and UMPA,
\'Ecole Normale Sup\'erieure de Lyon, 46 All\'ee d'Italie, 69364 Lyon, France}
\email{abdelghani.zeghib@ens-lyon.fr}

\date{\today}

\begin{abstract}
We present a list of open questions on various aspects of AdS geometry, that is, the geometry
of Lorentzian spaces of constant curvature $-1$. When possible we point out relations
with discrete subgroups of Lie groups and their actions on homogeneous spaces, with Teichm\"uller theory, as well as 
analogs in hyperbolic geometry.
\end{abstract}

\maketitle

\tableofcontents

\section{Introduction}

Let $n\in \N, n\geq 2$. The $n$-dimensional \emph{anti-de Sitter} (AdS) space, $\AdS_n$, is a complete Lorentzian
space of constant sectional curvature $-1$. To define it, let $\R^{n-1,2}$ be the vector space $\R^{n+1}$ endowed with
a symmetric bilinear form of signature $(n-1,2)$:
$$ b(x,y) = - x_1y_1 - x_2 y_2 + x_3y_3 + \cdots + x_{n+1}y_{n+1}~. $$
Then $\AdS_n$ can be defined as the following quadric in $\R^{n-1,2}$, with the Lorentzian metric induced by~$b$:
$$\AdS_n = \{ x\in \R^{n-1,2}~|~ b(x,x)=-1\}\,. $$
The isometry group of $\AdS_n$ is $\OO(n-1,2)$, and $\AdS_n$ identifies with $\OO(n-1,2)/\OO(n-1,1)$.
It is actually common to define $\AdS_n$ rather as the quotient of this quadric by the antipodal map $x\mapsto -x$ (corresponding to $-1\in\OO(n-1,2)$), so that we have the projective realization
$$\AdS_n = \{ [x]\in\R P^n : b(x,x)<0\}\,,$$
with boundary at infinity $\partial_{\infty}\AdS_n=\{ [x]\in\R P^n : b(x,x)=0\}$.
There is a formal analogy between $\AdS_n$ and the $n$-dimensional hyperbolic space $\HH^n$, and $\AdS_n$ can
indeed be considered, in several ways, as a Lorentzian analog of $\HH^n$: for instance, $\AdS_n$ and $\HH^n$ have the same complexification (see Section~\ref{sec:complex}).
Some of the questions below make sense
in this perspective. In particular, Sections \ref{sec:convex} and~\ref{sec:volume} are based on a strong similarity between quasi-Fuchsian hyperbolic 3-manifolds and globally hyperbolic AdS 3-manifolds.

However $\AdS_n$ was first introduced by physicists studying gravitation, as a possible model of the
universe --- a solution of the Einstein equations with negative cosmological constant. This physical
point of view provides another motivation for AdS geometry. From this perspective it seems natural
to study globally hyperbolic spaces, but also other geometric objects motivated by physics, such as
multi-black holes, ``particles'', or time functions. This point of view is apparent in Sections \ref{sec:time}, \ref{sec:blackholes}, and~\ref{sec:tachyons}. 

A third motivation comes from the theory of quotients (or Clifford--Klein forms) of homogeneous spaces, in particular affine symmetric spaces, by discrete groups acting properly discontinuously.
The case of Riemannian locally symmetric spaces is well understood, but much less is known about quotients of homogeneous spaces admitting a pseudo-Riemannian but no Riemannian structure.
It is natural to start by looking at the Lorentzian case, in particular in constant curvature.
Section~\ref{sec:compact} below mostly follows this direction.

A fourth motivation comes from relations between AdS geometry and Teichm\"uller theory. Mess
\cite{mess} discovered a simple proof of Thurston's Earthquake Theorem based on geometric
properties of globally hyperbolic AdS 3-manifolds, and, more recently, several 
results in Teichm\"uller theory have been obtained using arguments based on AdS geometry. 
This is a motivation for many of the questions stated in Sections \ref{sec:convex}, \ref{sec:volume}, and~\ref{sec:max}. 

We do not intend to state any new result here, but rather to present together a wide
spectrum of current open questions on AdS geometry. 



\section{Structure of compact AdS manifolds}\label{sec:compact}

In this section we consider AdS manifolds that are compact without boundary.
Although compact spacetimes are not physically relevant, they may lead to interesting mathematical questions. 
We first recall some key questions in the more general context of locally homogeneous manifolds.

\subsection{General questions on locally homogeneous manifolds}

Compact AdS manifolds are some of the simplest examples of non-Riemannian locally homogeneous manifolds, 
that is, of manifolds $M$ that are locally modeled on a homogeneous space $X=G/H$ 
where the left $G$-action does not preserve a Riemannian metric (this implies that $H$ is not compact).
Here are three important questions on such manifolds; the second one is more geometric, the third one more algebraic.\\

\noindent
{\bf Existence:} Let $X=G/H$ be a non-Riemannian homogeneous space.
Does there exist a compact manifold~$M$ locally modeled on~$X$? \\

\noindent
{\bf Completeness:} If $M$ is such a manifold, is it necessarily covered by a universal covering of~$X$?
In the case when $X$ is simply connected, we are asking whether $M$ necessarily identifies with the quotient $\Gamma\backslash X$ of~$X$ by some discrete subgroup~$\Gamma$ of~$G$ acting properly discontinuously.
This is equivalent to a notion of geodesic completeness when the $G$-action preserves a pseudo-Riemannian metric. \\

\noindent
{\bf Standardness (or $\pi_1$-rigidity):} If $M$ is a complete manifold locally modeled on~$X$, does it identify (up to finite covering) with a quotient $\Gamma\backslash X$ where $\Gamma$ is contained in some Lie subgroup~$L$ of~$G$ acting properly on~$X$?
(In this case, $L$ necessarily acts cocompactly on~$X$ and $\Gamma$ is a uniform lattice in~$L$.)
Such a quotient is called \emph{standard}.\\

We note that since $H$ is non-compact, the full group~$G$ does not act properly on~$X$; 
discrete subgroups~$\Gamma$ of~$G$ acting properly discontinuously cannot be uniform lattices in~$G$.
In the case of unimodular affine structures, \emph{i.e.}\ for $G=\SL_n(\R)\ltimes\R^n$ and $X=\R^n$, 
the completeness question is known as the Markus conjecture and standardness is related to the Auslander conjecture, 
which generalizes Bieberbach's theory of crystallographic groups (see \cite{abe01}).

\subsection{Nonstandard compact AdS manifolds} \label{ssc:nonstandard}

We now consider the previous questions for the anti-de Sitter space $X=\AdS_n=\OO(n-1,2)/\OO(n-1,1)$.\\

\noindent
{\bf Existence:} Compact AdS manifolds do not exist in even dimension, 
because a generalized version of the Gauss--Bonnet formula would imply that the Euler characteristic of such manifolds is non-zero, 
which would contradict the Poincar\'e--Hopf theorem.
However, compact AdS manifolds exist in any odd dimension $n=2m+1$, as was first observed by Kulkarni \cite{kul81}.
Indeed, by realizing $\AdS_{2m+1}$ as a real quadric as in the introduction and then identifying $\R^2$ with~$\C$, we see that the group $L=\U(m,1)$ acts transitively on $\AdS_{2m+1}$, with compact stabilizers, and that this action preserves the Lorentzian structure; therefore, any torsion-free uniform lattice~$\Gamma$ of $\U(m,1)$ yields a standard compact manifold $\Gamma\backslash\AdS_{2m+1}$ (such lattices exist by classical results of Borel and Harish-Chandra \cite{bor63}).\\

\noindent
{\bf Completeness:} 
All compact AdS manifolds are complete: this was proved by Klingler \cite{kli96}, following work of Carri\`ere \cite{car89} on the flat case.
We note that $\AdS_n$ is not simply connected: it is homeomorphic to $\R^{n-1}\times\mathbb{S}^1$.
Kulkarni--Raymond \cite{kr85} and Zeghib \cite{zeg98}~proved that compact AdS manifolds all have \emph{finite level}, which means that they are quotients $\Gamma\backslash\AdS_n$ of $\AdS_n$ up to finite covering (it is not necessary to consider a universal covering of $\AdS_n$).\\

\noindent
{\bf Standardness:}
We just saw that standard compact AdS manifolds exist in any odd dimension $n=2m+1$: we may take $L=\U(m,1)$.
There is a discrepancy between dimension~$3$ and higher dimension.
In dimension~$3$, non-standard compact AdS manifolds also exist: 
the first examples were constructed by Goldman \cite{Goldman.Non.Standard} by deformation of standard manifolds, 
and a systematic study was then carried out by Salein \cite{sal00}, Kassel \cite{kasselPhD}, and Gu\'eritaud--Kassel \cite{gueritaud-kassel}.
The existence of non-standard compact manifolds and the possibility to deform all compact AdS manifolds 
in dimension~3 are linked to the fact that the isometry group $\OO(2,2)$ of $\AdS_3$ is not simple: it is the direct product $\PSL_2(\R)\times\PSL_2(\R)$ up to finite index and a $2$-fold covering.
On the other hand, the isometry group $\OO(n-1,2)$ of $\AdS_n$ is simple for $n>3$, 
and the answer to the following question (first asked by Zeghib \cite{zeg98}) is conjectured to be affirmative.

\begin{question}\label{standardcompact}
Suppose $n=2m+1>3$.
Let $M=\Gamma\backslash\AdS_n$ be a compact AdS manifold, 
where $\Gamma$ is a discrete group of isometries acting properly discontinuously on $\AdS_n$.
Is $\Gamma$ necessarily a uniform lattice in $\U(m,1)$ (up to conjugacy)?
\end{question}

Question~\ref{standardcompact} asks not only whether the proper $\Gamma$-action extends to a proper action of a Lie group~$L$ (\emph{i.e.}\ whether $M$ is standard), but also whether this $L$ is in fact conjugate to $\U(1,m)$.
By \cite{zeg98}, if not, then $M$ is non-standard and $\Gamma$ is Zariski-dense in $\mathrm{(S)O}(n-1,2)$.
By cohomological arguments of Raghunathan \cite{rag65} and Weil \cite{wei64}, 
standard compact AdS manifolds in dimension $n>3$ cannot be deformed into non-standard ones.

Here is a weaker question than Question~\ref{standardcompact}, stating the uniqueness of AdS structures on standard compact manifolds in dimension $>3$: is it true that for $m>1$, no uniform lattice of $\U(m,1)$ can be embedded into $\OO(2m,2)$ as a Zariski-dense subgroup acting properly discontinuously on $\AdS_{2m+1}$?
A positive answer to this question would be implied by a positive answer to Question~\ref{standardcompact} because of the following fact \cite{ser71}: if $\Gamma$ is a discrete group of isometries acting properly discontinously on $\AdS_{2m+1}$, then the quotient $\Gamma\backslash\AdS_{2m+1}$ is compact if and only if the cohomological dimension of a torsion-free finite-index subgroup of~$\Gamma$ is $2m$ (which is true when $\Gamma$ is a uniform lattice in $\U(m,1)$).
This is a case where the Margulis superrigidity theorem on representations of lattices of semisimple Lie groups does not apply.

\subsection{Moduli space and volume of compact AdS manifolds in dimension~3}\label{subsec:compactAdS3}

We now consider dimension~3.
As we just explained, compact AdS 3-manifolds are all quotients 
$\Gamma\backslash\AdS_3$ of $\AdS_3$ by a discrete group of isometries~$\Gamma$ acting properly discontinuously, up to finite cover; the isometry group $\OO(2,2)$ of $\AdS_3$ is $\PSL_2(\R)\times\PSL_2(\R)$ up to finite index and a $2$-fold covering.
By Kulkarni--Raymond \cite{kr85}, all torsion-free discrete subgroups~$\Gamma$ of $\PSL_2(\R)\times\PSL_2(\R)$ acting properly discontinuously and cocompactly on $\AdS_3$ are of the form
$$\Gamma = \{ (j(\gamma),\rho(\gamma)) : \gamma\in\pi_1(S)\} ,$$
where $S$ is a closed hyperbolic surface and $j$ and~$\rho$ are representations of $\pi_1(S)$ into $\PSL_2(\R)$ with one of them Fuchsian (\emph{i.e.}\ injective and discrete).
However, not all pairs $(j,\rho)\in\nolinebreak\Hom(\pi_1(S),\PSL_2(\R))^2$ with $j$ Fuchsian are \emph{admissible}, \emph{i.e.}\ lead to properly discontinuous actions on~$\AdS_3$: by \cite{kasselPhD}, a necessary and sufficient condition is that
$$C_{length}(j,\rho) := \sup_{\gamma\in\pi_1(S)\smallsetminus\{ e\}} \frac{\lambda(\rho(\gamma))}{\lambda(j(\gamma))} < 1,$$
or equivalently that
$$C_{Lip}(j,\rho) := \inf\big\{ \mathrm{Lip}(f)\ :\ f : \HH^2\rightarrow\HH^2\text{ $(j,\rho)$-equivariant}\big\} < 1,$$
where $\mathrm{Lip}$ denotes the Lipschitz constant and for $g\in\PSL_2(\R)$ we set $\lambda(g):=\inf_{x\in\HH^2} d_{\HH^2}(x,g\cdot x)$ (this is the translation length of~$g$ if $g$ is hyperbolic, and $0$ if $g$ is parabolic or elliptic).
In particular, if $(j,\rho)$ is admissible and $j$ Fuchsian, then $\rho$ cannot be Fuchsian.
Indeed, Thurston \cite{thu86} proved that $C_{length}$ and $C_{Lip}$ are non-negative on $\cT\times\cT$, where $\cT$ denotes the Teichm\"uller space of~$S$ (this is the starting point of the theory of the \emph{Lipschitz metric} on~$\cT$ developed in \cite{thu86}).

Understanding the moduli space of compact AdS manifolds reduces to understanding, for a fixed closed hyperbolic surface~$S$, the space~$\mathrm{Adm(S)}$ of admissible pairs $(j,\rho)\in\Hom(\pi_1(S),\PSL_2(\R))^2$ with $j$ Fuchsian.
It follows from completeness and from the Ehresmann--Thurston principle that $\mathrm{Adm(S)}$ is open in $\Hom(\pi_1(S),\PSL_2(\R))^2$.
Moreover, Salein \cite{sal00} proved (by some explicit geometric construction) that there exist admissible pairs $(j,\rho)\in\mathrm{Adm}(S)$ such that $\rho$ may belong to any connected component of $\Hom(\pi_1(S),\PSL_2(\R))$ except for the two extremal ones (corresponding to Fuchsian representations).
Thus $\mathrm{Adm}(S)$ has at least $4g-5$ connected components.

\begin{question}
How many connected components does the space $\mathrm{Adm}(S)$ have?
What is its topology?
\end{question}

Another question concerns the volume of compact AdS manifolds.

\begin{question}
Which values can the volume of a compact AdS $3$-manifold take?\\
For $\Gamma=\{ (j(\gamma),\rho(\gamma)) : \gamma\in\pi_1(S)\}$ with $(j,\rho)\in\mathrm{Adm}(S)$, does the volume of $\Gamma\backslash\AdS_3$ depend only on the Euler number of~$\rho$, or does it vary when $(j,\rho)$ varies in a given connected component of $\mathrm{Adm}(S)$?
\end{question}

\section{Convex AdS manifolds with boundary in dimension 3}\label{sec:convex}

We now turn our attention to \emph{globally hyperbolic} AdS 3-manifolds. Recall that a Lorentzian manifold
is said to be globally hyperbolic if it contains a \emph{Cauchy hypersurface}, 
\emph{i.e.}\ a space-like hypersurface $S$ such that all inextendible time-like
lines intersect $S$ exactly once. It is \emph{globally hyperbolic maximal compact} (GHMC) if in addition $S$ is closed
and any isometric embedding of $M$ into another Lorentzian manifold of the same dimension with the same property 
is an isometry.

In this section, we denote by $S$ a closed surface with genus $g\geq 2$ and by $\cT$ the Teichm\"uller space of $S$.
We set $M=S\times\R$ and denote by $\cGH$ the space of GHMC AdS structures on $M$, considered up to isotopy.

The questions below can be related to two classical geometric results.
Firstly, Alexandrov \cite{alex} proved that if $h$ is a hyperbolic metric on the sphere $\mathbb{S}^2$ 
with cone singularities of angle less than $2\pi$, then $h$ defines the same length space distance on $\mathbb{S}^2$ as
a unique polyhedral embedding sending $\mathbb{S}^2$ to the the boundary of a convex polyhedron in $\HH^3$ --- in simpler terms,
$h$ is realized on the boundary of a unique convex polyhedron in $\HH^3$. 
In the context of AdS geometry it is quite natural to replace
the convex polyhedra homeomorphic to balls considered by Alexandrov by AdS spaces with more topology and a 
space-like boundary --- this leads to the questions in Section~\ref{ssc:metric} and below.
Secondly, Andreev \cite{andreev,andreev-ideal}
characterized the possible dihedral angles of convex polyhedra in $\HH^3$, and, in the AdS context, this also leads
quite naturally to the questions considered in Section \ref{ssc:lamination} and to other questions mentioned
below on the angles of measured laminations of AdS manifolds with space-like, convex boundary.

\subsection{Prescribing the measured bending lamination on the boundary of an AdS convex core}
\label{ssc:lamination}

Let $h\in \cGH$ be a GHMC AdS structure on $M$. Then $(M,h)$ contains a unique smallest non-empty 
convex subset, its \emph{convex core}, which we shall denote by $C(M)$.
Here we say that a subset $K\subset M$ is convex if any geodesic segment of $M$ with endpoints in $K$ is contained in $K$.
The boundary of $C(M)$ is the union
of two bent surfaces $S_-, S_+$, each homeomorphic to $S$, 
with measured bending laminations $l_-, l_+$, see \cite{mess,mess-notes}. 
(Note that $S_-$ and $S_+$ are disjoint unless $h$ is Fuchsian, in which case $S_-=S_+$ is a totally 
geodesic surface in $M$.) Moreover, it is not difficult to prove that $l_-$ and $l_+$ fill $S$.

\begin{question} \label{q:lamination}
If $l_-$ and $l_+$ are two measured laminations which fill $S$, is there a unique GHMC AdS structure
$h$ on $M$ such that $l_-$ and $l_+$ are the measured bending laminations of the boundary of the 
convex core?
\end{question}

The existence of $h$ is proved in \cite{earthquakes}, however the uniqueness remains open. This is
quite similar to the corresponding question on convex cores of geometrically finite hyperbolic
3-manifolds, for which the existence is proved in \cite{bonahon-otal} but the uniqueness is
still open.

Question \ref{q:lamination} can be translated purely in terms of hyperbolic metrics and
earthquakes, as pointed out in \cite{earthquakes}. The following is an equivalent formulation.

\begin{question} \label{q:fixedpoint}
Let $l_-, l_+$ be two measured laminations which fill $S$. Does the composition 
$E_l(l_-)\circ\nolinebreak E_l(l_+)$ have a unique fixed point in $\cT$?
\end{question}

Here $E_l(l_{\pm})$ is the right earthquake along $l_{\pm}$, considered as a homeomorphism of $\cT$.
Again the existence of a fixed point is proved in \cite{earthquakes}, but not its uniqueness.

\subsection{Prescribing the induced metric on the boundary of an AdS convex core}
\label{ssc:metric}

Question~\ref{q:lamination} has an analog concerning the induced hyperbolic metrics $m_-,m_+\in\cT$ on the two connected
components of the boundary of the convex core.

\begin{question} \label{q:metrics}
Let $m_-, m_+\in \cT$. Is there a unique $h\in \cGH$ such that the induced metric on the two
connected components of the boundary of the convex core are $m_-$ and $m_+$?
\end{question}

The existence part is (almost) known, it is obtained in the PhD thesis of Boubacar 
Diallo\footnote{To appear soon.}. The uniqueness however remains open. 

This is again analogous to the situation of geometrically finite hyperbolic 3-manifolds,
where the existence is known --- it can be obtained as a consequence of results of \cite{L4},
or of \cite{epstein-marden}. 

Question \ref{q:metrics} can also be translated into a purely 2-dimensional question. Given two
hyperbolic metrics $h_l, h_r\in \cT$, there is a unique left earthquake and a unique right
earthquake sending $h_l$ to $h_r$. Together, they form an earthquake ``bigon'' with vertices 
$h_l$ and $h_r$, and we consider the midpoints $h_+$ and $h_-$ of the two edges. Let 
$U:\cT\times \cT\rightarrow \cT\times\cT$ be the map sending $(h_l,h_r)$ to 
$(h_-, h_+)$. 

\begin{question}
Is $U$ a homeomorphism? 
\end{question}

This question is actually equivalent to Question \ref{q:metrics}.

Note that the restriction of $U$ to the diagonal is the identity, however it is not even
clear whether $U$ is a homeomorphism in a neighborhood of the diagonal. Investigating this
question leads to a simple infinitesimal statement which does not seem obvious.

\subsection{The boundary of the convex core with particles}

Consider now a globally hyperbolic AdS 3-manifold containing massive particles --- that is,
cone singularities along time-like lines --- of angle less than $\pi$. 
It is shown in \cite{cone} that these manifolds have a well-defined convex core. Moreover,
the particles do not intersect the support of the measured bending lamination of the 
boundary of the convex core, and they are orthogonal to the boundary of the convex core. 

It follows that the induced metric on the boundary of the convex core is again hyperbolic,
with cone singularities at the intersections with the particles, of angle equal to the
angle of the 3-manifold at those particles. 

One can therefore ask both Question \ref{q:lamination} and Question \ref{q:metrics} in this
context (again, as long as the angles are less than $\pi$). Note that some of the arguments
that could possibly be used to answer Questions \ref{q:lamination} and \ref{q:metrics} 
will not be available when particles are present.

\subsection{AdS manifolds with smooth, locally convex boundary}

Consider now a convex subset $K$ of a GHMC AdS 3-manifold $M$. Suppose that
$K$ has a smooth, strictly convex boundary --- by strictly convex we mean here that
the second fundamental form of the boundary is positive definite at each point.
The Gauss formula then shows that the induced metric on the boundary of $K$
has curvature $\kappa<-1$ everywhere. Rather than as a subset, $K$ can be considered as an
AdS manifold with smooth, space-like, locally strictly convex boundary.

\begin{question} \label{q:I-smooth}
Let $h_-, h_+$ be two metrics on $S$ with curvature $\kappa<-1$. Is there a unique
AdS manifold $N$ with smooth, space-like, strictly convex boundary such that the induced
metrics on the two connected components of $\partial N$ are equal respectively to $h_-$ and $h_+$?
\end{question}

Note that the existence part of this question might be much easier than the
uniqueness. It is conceivable that the type of argument used in \cite{L4} could
be used in the AdS setting considered here.

The analogous question for hyperbolic manifolds has a positive answer \cite{hmcb}.
A positive answer is also known \cite{iie} when $h_-=h_+$, corresponding to Fuchsian manifolds,
that is, AdS manifolds containing a closed, totally geodesic space-like surface.

Using the duality between convex space-like surfaces in $\AdS_3$, one finds an 
equivalent form of this question where one wants to prescribe the third fundamental
form, rather than the induced metric, on the boundary. This equivalent question is strongly 
related to Question \ref{q:lamination} above. In the hyperbolic case the
analogous statement is not equivalent to the statement on the induced metric,
but it holds, too \cite{hmcb}.

\subsection{AdS manifolds with polyhedral boundary}

The questions in the previous paragraph can also be asked in the polyhedral 
category, that is, when one considers AdS manifolds with convex, polyhedral
boundary. One added twist is that, to get a complete statement, one cannot
exclude the possibility that the boundary of such an AdS manifold intersects the boundary of its convex core; in other
terms, the boundary should not be restricted to be polyhedral, but should be
allowed to have some regions pleated along a measured lamination. 
If we understand the term ``AdS manifold with polyhedral boundary'' in this
manner, we can ask the following question.

\begin{question} \label{q:I-poly}
Let $h_-, h_+$ be two hyperbolic metrics on $S$ with cone singularities of
angle larger than $2\pi$. Is there a unique AdS manifold with convex,
polyhedral boundary for which the induced metric on the boundary is given
by $h_-$ and $h_+$?
\end{question}

Note that Question \ref{q:metrics} is a special case of Question~\ref{q:I-poly}, corresponding
to metrics with no cone singularity. Here, too, a positive answer has been obtained
by Fillastre \cite{fillastre3} in the Fuchsian case, when $h_-=h_+$.

Note also that it should be possible to consider both Question \ref{q:I-smooth} and
Question \ref{q:I-poly} as a special case of a question on the induced metrics on the
boundaries of AdS manifolds with convex boundary, with no hypothesis on the smoothness
nor the strict convexity of the boundary. A natural statement is that one could obtain
in this way any pair of $CAT(-1)$ metrics in exactly one way. The analogous 
statement might also hold for hyperbolic manifolds with convex boundary (with no
smoothness assumption); the metrics obtained should then be any metrics with
curvature larger than $-1$ in the sense of Alexandrov.

Question \ref{q:I-poly} has relations with the flippable tilings of hyperbolic
surfaces considered in \cite{fs}.

\section{Volume and width of the convex core in dimension 3}\label{sec:volume}

As in the previous section, we denote by $S$ a closed surface with genus $g\geq 2$ and by $\cT$ the Teichm\"uller space of $S$.
We set $M=S\times\R$ and denote by $\cGH$ the space of globally hyperbolic maximal compact AdS structures on $M$, considered up to isotopy.

\subsection{The volume of the convex core of AdS manifolds}

Given a quasi-Fuchsian hyperbolic metric on~$M$, one can consider two different volumes: 
its renormalized volume (see \cite{volume,review}) and the volume of the convex
core. Brock \cite{brock:2003} proved that the volume of the convex core is quasi-equivalent
to the Weil--Petersson distance between the induced metrics on the boundary
of the convex core. Moreover, the renormalized volume and the volume of the
convex core are quite close \cite{compare}, so that the renormalized volume
is also quasi-equivalent to the Weil--Petersson distance between the conformal
metrics at infinity.

For a globally hyperbolic AdS manifold $M$ one can consider the volume $V(C(M))$ of
the convex core $C(M)$, or the total volume $V(M)$ of $M$, which is finite. The two quantities
are actually related by a simple relation:
$$ V(C(M)) + V(M) = 2\pi^2(g-1) + \frac{\mathrm{length}_m(l)}{2}~, $$
where $m$ is the induced metric and $l$ the measured bending lamination of the boundary of the convex core.

Note that in the AdS setting $\mathrm{length}_m(l)$ is not bounded, so both $V(M)$ and
$V(C(M))$ are presumably unbounded, and they might be quite different.

\begin{question}
Is $V(M)$ or $V(C(M))$ quasi-equivalent to a distance (for instance the Teichm\"uller
distance, or the Weil--Petersson distance) between the induced metrics on the boundary of
the convex core, or between the left and right hyperbolic metrics? 
\end{question}

\subsection{Minima and convexity of the volume}

Note that for Fuchsian GHMC AdS manifolds, the volume $V(M)$ is equal to $\pi/2$ times the area of
any hyperbolic metric on $S$, that is, $2\pi^2(g-1)$. 

\begin{question}
Is $V(M)$ minimal at the Fuchsian locus?
\end{question}

More specifically, fix a hyperbolic metric $h_l\in \cT$ on $S$. For each $h_r\in \cT$
there is (see \cite{mess,mess-notes}) a unique GHMC AdS manifold $M(h_l, h_r)$ with left and right
hyperbolic metrics equal to $h_l$ and $h_r$, respectively, and we can consider the functions 
$V(M(h_l,\bullet))$ and $V(C(M(h_l, \bullet)))$ on~$\cT$. 

\begin{question}
Are these functions convex? Are they proper? Do they have a unique minimum at~$h_l$?
\end{question}

It would also be interesting to understand the ratio $V(C(M))/V(M)$ as a function either on 
$\cT\times \cT$, or on $\cT$ for $h_l$ fixed.

\subsection{The width of the convex core}

Another interesting quantity for a GHMC AdS manifold~$M$ is the width $w(M)$ of its convex core, \emph{i.e.}\ the maximum (time-like) distance of points of the convex core to the lower boundary (this
notion is used in \cite{maximal}). This width is always less than $\pi/2$. It would be
interesting to understand whether there is a good relation between this width and
the Teichm\"uller (or Weil--Petersson) distance between the left and right hyperbolic
metrics of~$M$. It follows from estimates in \cite{maximal} that, if the Teichm\"uller
distance between $h_l$ and $h_r$ is very large, then $w(M)$ is close to $\pi/2$.

There is a related question which can be stated purely in 2-dimensional terms. 
Let $h_l,h_r$ be two hyperbolic metrics on $S$. There exists a unique minimal Lagrangian
diffeomorphism between $(S,h_l)$ and $(S, h_r)$ which is isotopic to the identity (see \cite{schoen:role,L6}); let
$ml(h_l,h_r)$ be its distorsion. Then $ml(h_l, h_r)$ can be expressed in terms of the
supremum of the principal curvatures of the (unique) maximal space-like surface in the globally hyperbolic AdS manifold $M$ with left and right hyperbolic metrics equal to $h_l$ and $h_r$. 

Clearly, $ml(h_l, h_r)$ is larger than the Teichm\"uller distance between $h_l$ and $h_r$.
A natural question is whether a bound on the Teichm\"uller distance implies a bound on $ml(h_l,hr)$. 
In other words, do minimal lagrangian diffeomorphisms give an efficient way to measure the Teichm\"uller 
distance?

It appears possible that this question is related to the notion of width of $M$.

\section{Fundamental groups of globally hyperbolic AdS manifolds in higher dimension}

Another yet unsolved question about GHMC AdS manifolds concerns their possible topology, \emph{i.e.}\ the possible topology of their Cauchy (hyper)surfaces.

In the case of \emph{flat} GHMC manifolds, the problem has been solved: Scannell \cite{sca01} proved that in dimension $3+1$, every Cauchy hypersurface must be diffeomorphic to a compact quotient (by a discrete group of isometries) of, either the Euclidean $3$-space, the hyperbolic space $\HH^3$, or the product $\R\times\HH^2$. This result has been generalized to any dimension in \cite{bar05}, where it is proved that Cauchy hypersurfaces of GHMC AdS manifolds of dimension $n+1$ are diffeomorphic to a compact quotient (by discrete groups of isometries) of $\R^p\times\HH^q$ with $p+q=n$. 

In the case of GHMC manifolds locally modeled on the \emph{de Sitter} space $\mathrm{dS}_n$ (which is the analog of $\AdS_n$ in constant positive curvature), it follows from Scannell's thesis \cite{scaPhD} that the problem is equivalent to the problem of characterizing closed manifolds admitting conformally flat Riemannian metrics.

For the case under interest here, the natural result to be expected is the following.

\begin{question}\label{CauchyGHMC}
Prove that the Cauchy surfaces of a GHMC manifold locally modeled on $\AdS_{n+1}$ are diffeomorphic to the quotient of $\HH^p\times\HH^q$ by a uniform lattice of $\SO(1,p)\times\SO(1,q)$ with $p+q = n$.
\end{question}

It is actually easy to construct, given a uniform lattice $\Gamma$ of $\SO(1,p)\times\SO(1,q)$, a globally hyperbolic AdS spacetime with Cauchy hypersurfaces diffeomorphic to $\Gamma\backslash (\HH^p\times\HH^q)$.
Indeed, $\Gamma$ can be embedded into $\SO(2,n)$ as a discrete subgroup through the natural embedding $\SO(1,p)\times\nolinebreak\SO(1,q)\hookrightarrow\SO(2,n)$; if $\Omega$ denotes the domain of $\AdS_{n+1}$ whose points are not causally related to any $\Gamma$-iterate of themselves, then the quotient $\Gamma\backslash\Omega$ is a GHMC manifold with the required property. We call these examples \emph{$GH$-standard}.
Question~\ref{CauchyGHMC} would be a corollary of the following more difficult question.

\begin{question}
Is every GHMC AdS manifold a deformation of a $GH$-standard AdS manifold?
\end{question}

We note that certain uniform lattices of $\SO(1,n)$ do admit non-trivial (and even Zariski-dense) deformations into $\SO(2,n)$, by a bending construction due to Johnson and Millson \cite{jm87} (see \cite{kas12} for details).
However, it is a non-trivial problem to characterize the uniform lattices of $\SO(1,n)$ that admit such non-trivial deformations; this problem is similar to the maybe better-known problem of their deformation into $\SO(1,n+1)$ (see \cite{jm87,apa92,tan93}).

On the other hand, it is not clear whether some uniform lattices of $\SO(1,p)\times\SO(1,q)$ 
for $p,q\geq 2$ admit non-trivial deformations into $\SO(2,n)$.
Such lattices would not be irreducible, by the Margulis superrigidity theorem.
For $p=1$, the group $\SO(1,p)$ is isomorphic to~$\R$, which allows for a rich deformation theory.

\section{Geometric time functions on AdS manifolds}\label{sec:time}

The simplest examples of GHMC spacetimes are metric products, \emph{i.e.}\ spacetimes of the form $M= (\R,-dt^2)\oplus\nolinebreak (\Sigma, h)$ where $(\Sigma,h)$ is a compact Riemannian manifold.  
An interesting generalization is those where the metric has a warped form $-dt^2 + w(t) h$.
In fact, from a topological viewpoint, any GHMC spacetime $M$ is homeomorphic to a product $\R\times\Sigma$, where $\Sigma$ is a Cauchy hypersurface of $M$, but in general $M$ is far from being a product from a metric point of view.
In order to evaluate the metric distortion, it is natural and worthwhile to ask if there are privileged splittings 
(or \emph{time functions}) for a given GHMC spacetime. 
It is especially exciting to ask what  geometry remains from the couple of the orthogonal foliations after perturbation 
of a direct or a warped  product metric $-dt^2\oplus w(t)h$.
So the general natural question one can ask is to produce canonical geometric foliations by Cauchy hypersurfaces, 
or equivalently geometric time functions, which yield a kind of measurement of the default for $(M,g)$ to be a metric (warped) product.  
Actually, asymptotic behavior, singularities, shocks, and similar  questions are meaningful only in a ``natural'' coordinates  
system (which does not create artificial singularities).

\subsection{$F$-time functions}

It turns out that a reasonable natural requirement for special time functions is to be solutions of a geometric PDE, as in the case of the so-called $F$-time functions.
Here $F$ stands for a function
$F: {\R}^n \to {\R}$  which is invariant by permutation of the coordinates. 
The $F$-curvature of a hypersurface $\Sigma$ embedded in an $(n+1)$-dimensional spacetime $M$ is defined 
by evaluating the function $F$ on the eigenvalues of the second fundamental form of $\Sigma$.
An $F$-time function is a time function $t:M  \to \R$ such that each level set $t^{-1}(a)$ has constant $F$-curvature, and this curvature is increasing
with $a$.

A \emph{constant mean curvature} (CMC) time corresponds to the  case where $F$ is the (arithmetic) mean  of eigenvalues, and \emph{K-times} are defined by taking $F$ to be the opposite of the product of the eigenvalues;  
the $K$-curvature is sometimes called the Gauss--Killing--Kronecker--Lipschitz curvature.

One interesting feature is that if an $F$-time exists then it is unique and any Cauchy $F$-hypersurface coincides with
a level set of this time. Furthermore, in general,  existence of such a time allows one to solve the problem of prescribing the $F$-curvature.

Observe that among GHMC AdS manifolds, only the Fuchsian ones  have a  warped product structure.
The associated time is in fact their cosmological time and it is an $F$-time for any~$F$.
Outside the Fuchsian case, all theses times are different, if they exist.  
This motivates the question of their existence as well as their asymptotic comparison (at least for perturbation of Fuchsian manifolds).

For GHMC AdS manifolds (and their cousins of constant non-negative sectional curvature), existence of CMC times was proved in \cite{BBZ,ABBZ}.

As for $K$-times, the question of their  existence was handled in \cite{BBZ2} but only in dimension 3. 
In fact, GHMC AdS manifolds do not support  (everywhere-defined) $K$-times. 
It is the complement of the convex core which admits such a time. 
More precisely, the past of the convex core admits a $K$-time, with a $K$-curvature of leaves having range $(-\infty,0)$. 
(Its future has a reversed $K$-time). In general dimension, one can ask the following question. 

\begin{question} 
Does the past of the convex core of a GHMC AdS manifold (of arbitrary dimension) always possess a $K$-time?
\end{question}

The question may be posed for more general $F$-curvatures, \emph{e.g.}\ when
$F: \R^n \to \R$ is a symmetric polynomial. For instance, the quadratic case $F (\lambda_1, \ldots \lambda_n) =
\Sigma \lambda_i \lambda_j$ corresponds, up to  constants, to the scalar curvature.

In the similar flat situation, $K$-times were used in \cite{BBZ2} to solve an equivariant Minkowski problem
in the 3-dimensional Minkowski spacetime. Higher-dimensional variants were in particular considered in
\cite{Bertrand, Fillastre}.
 
\subsection{Special foliations with cone singularities}
 
Consider a globally hyperbolic AdS manifold $M$ with particles of cone angle less than $\pi$, as
seen in \cite{cone}. 

\begin{question} \label{q:maximal}
Is there a unique maximal space-like surface in $M$ orthogonal to the particles?  
\end{question}

This would extend a result of \cite{minsurf}. A positive answer would probably lead to a
positive answer to the following, strongly related question (following the arguments
in \cite{minsurf}).

\begin{question}
Let $S$ be a closed surface.
Choose $x_1,\cdots, x_n\in S$ and $\theta_1,\cdots, \theta_n\in (0,\pi)$, and let $h,h'$ be two hyperbolic metrics on $S$ with cone singularities of angle $\theta_i$ at $x_i$. Is there a unique minimal Lagrangian diffeomorphism isotopic to the identity fixing the $x_i$ between $(S,h)$ and $(S, h')$? 
\end{question}

A natural extension of Question \ref{q:maximal} is whether this (supposedly unique) maximal surface
is the ``center'' leaf of a foliation of $M$ by constant mean curvature surfaces. 

\begin{question}
Does $M$ admit a unique CMC time with each surface of constant time orthogonal to the singularities?
\end{question}

This would be a natural extension, to manifolds with cone singularities, of the main result of 
\cite{BBZ}. Similarly, a positive answer to the following question would provide an extension of
the main result of \cite{BBZ2}.

\begin{question}
Does the past (resp.\ future) of the convex core of $M$ have a unique foliation by constant Gauss
curvature?
\end{question}

\section{Multi-black holes and complete surfaces}\label{sec:blackholes}

\subsection{Multi-black holes with particles}

Multi-black holes are 3-dimensional AdS manifolds which can be considered as analogs of globally
hyperbolic manifolds but based on a complete hyperbolic surfaces with infinite area, rather than 
on a closed hyperbolic surface, see \cite{Multi,barbot-1,barbot-2}. The space of multi-black
holes homeomorphic to $S\times \R$, where $S$ is the interior of a compact surface with 
non-empty boundary, is parameterized by the product of two copies of the Teichm\"uller space
of complete hyperbolic metrics on $S$.

It appears that many properties of globally hyperbolic AdS manifolds extend to the setting
of globally hyperbolic manifolds with ``particles'' --- cone singularities along time-like geodesics --- as
long as the cone angles are less than $\pi$, see \emph{e.g.}\ \cite{cone}. Note that this angle condition
prevents the singularities from ``interacting'', that is, the singular locus of these globally 
hyperbolic manifolds has to be a disjoint union of segments. The space of globally hyperbolic
metrics with $n$ particles of fixed angle on $S\times \R$ is then parameterized by the product of two copies
of $\cT_{S,n}$, the Teichm\"uller space of $S$ with $n$ marked points.

It might appear natural to try to extend this analysis to multi-black holes containing
particles, at least as long as the cone angles at the particles are less than $\pi$. 
One could imagine that the space of metrics of this type with $n$ particles of fixed angle 
on $S\times \R$, where $S$ is now the interior of a compact surface with non-empty boundary, 
is parameterized by the product of two copies of the Teichm\"uller space of $S$ with $n$
marked points. It might be relevant in this context to suppose that the particles enter
the compact part of the multi-black hole, or that they do not enter the asymptotic regions.

Using the argument developed by Mess \cite{mess}, it might then be possible to obtain as a 
corollary a version of Thurston's earthquake theorem for hyperbolic surfaces with 
boundary components and cone singularities of angle less than $\pi$, extending the
results of both \cite{cone} and \cite{mbh}.

A further step in this direction would be to extend to multi-black holes the (still
partial) results known for globally hyperbolic manifolds when the angles are only 
less than $2\pi$ and the particles are allowed to interact, see \cite{colI}.

\subsection{The convex core in the case of a proper action on $\AdS_3$}\label{subsec:coreproper}

We now consider the case without singularities.
Let $S$ be the interior of a compact surface with finitely many boundary components and let $\T_S$ be its Teichm\"uller space, \emph{i.e.}\ the set of (complete) convex cocompact hyperbolic structures on~$S$ up to isotopy.
Here ``convex cocompact'' means that we allow for funnels but not for cusps.
For $h_l,h_r\in\T_S$, let $\rho_l,\rho_r : \pi_1(S)\rightarrow\PSL_2(\R)$ be the corresponding holonomies.
By \cite{barbot-1,barbot-2,mbh}, the group $(\rho_l,\rho_r)(\pi_1(S))\subset\PSL_2(\R)\times\PSL_2(\R)$ acts properly discontinuously on some nonempty convex domain in $\AdS_3$, namely the interior of the convex hull of the graph in $\partial_{\infty}\AdS_3\simeq\R P^1\times\R P^1$ of a homeomorphism $\R P^1\rightarrow\R P^1$ conjugating $\rho_l$ to~$\rho_r$; Mess's theory (as in Sections \ref{sec:convex} and~\ref{sec:volume} above) extends to this context.
On the other hand, for certain pairs $(h_l,h_r)\in\T_S\times\T_S$ the action of $(\rho_l,\rho_r)(\pi_1(S))$ is properly discontinuous on \emph{the whole of $\AdS_3$}: by \cite{kasselPhD}, a necessary and sufficient condition is that
$$C_{length}(h_l,h_r) := \sup_{\gamma\in\pi_1(S)\smallsetminus\{ e\}} \frac{\mathrm{length}_{h_r}(\gamma)}{\mathrm{length}_{h_l}(\gamma)} < 1,$$
or equivalently that
$$C_{Lip}(h_l,h_r) := \inf_{\varphi\in\mathrm{Homeo}_0(S)} \mathrm{Lip}_{h_l,h_r}(\varphi) < 1$$
(up to exchanging $h_l$ and $h_r$), where $\mathrm{Lip}_{h_l,h_r}$ denotes the Lipschitz constant, measured in the initial metric~$h_l$ and the final metric~$h_r$.
Given $h_l\in\T_S$, it is always possible to construct $h_r\in\T_S$ such that these conditions are satisfied.
By \cite{gueritaud-kassel}, proper discontinuity on $\AdS_3$ is equivalent  to $C_{Lip}<1$ even if we allow for cusps.

\begin{question}\label{readproperness}
Let $h_l,h_r\in\T_S$ and let $\rho_l,\rho_r : \pi_1(S)\rightarrow\PSL_2(\R)$ be the corresponding holonomies.
How can one read the properness or non-properness of the action of $(\rho_l,\rho_r)(\pi_1(S))$ on $\AdS_3$:
\begin{itemize}
  \item on the hyperbolic metrics of the boundary components of the convex core?
  \item on the bending laminations?
\end{itemize}
\end{question}

\section{Maximal surfaces}\label{sec:max}

\subsection{Symplectic maps from maximal surfaces}

Let $M$ be a 3-dimensional globally hyperbolic AdS manifold. Following \cite{mess} we can associate
to it two points in Teichm\"uller space, corresponding to its left and right hyperbolic metrics
$h_l, h_r$. Mess \cite{mess} proved that this provides a parameterization of the space $\cGH$ of 
globally hyperbolic maximal compact AdS metrics on $M$ by $\cT\times \cT$, where $\cT$ is the Teichm\"uller
space of any closed space-like surface in $M$. 

Moreover, $M$ contains a unique closed, space-like maximal surface $S_0$. The second fundamental
form of $S_0$ is the real part of a holomorphic quadratic differential $q$ (for the complex structure
$c$ of the induced metric) and this defines a point $(c,q)\in T^*\cT$. This provides another
parameterization of the space $\cGH$ by $T^*\cT$, see \cite{minsurf}.

Composing these two parameterizations of $\cGH$, we obtain a map $\phi:T^*\cT\rightarrow \cT\times \cT$.

\begin{question}
Is $\phi$ symplectic, that is, is the pull-back by $\phi$ of the difference of the Weil--Petersson
symplectic forms on the two factors in the image equal (up to a constant) to the cotangent
symplectic structure on $T^*\cT$?  
\end{question}

This question has an extension to the universal Teichm\"uller space, in the setting considered in \cite{krasnov-scarinci}.

\subsection{Harmonic extensions of quasi-symmetric homeomorphisms of the circle}

R. Schoen \cite{schoen:role} asked the following question.

\begin{question}
Let $\psi : \mathbb{S}^1\rightarrow\mathbb{S}^1$ be an orientation-preserving quasi-symmetric homeomorphism. 
Is there a unique quasi-conformal harmonic diffeomorphism $\phi : \HH^2\rightarrow\HH^2$ such that $\phi$
extends continuously at infinity and its restriction to the boundary at infinity of 
$\HH^2$ is $\psi$? 
\end{question}

The uniqueness is known. The corresponding question for hyperbolic metrics on closed
surfaces has a positive answer (see \cite{schoen:role,L5}): given two hyperbolic
metrics $h$ and $h'$ on a closed surface $S$, there is a unique harmonic diffeomorphism
isotopic to the identity between $(S,h)$ and $(S,h')$. 

This question can be stated in terms of maximal surfaces in AdS, as should be clear from 
\cite{maximal} --- the interested reader should be able to obtain the equivalent formulation
without great difficulty.

\subsection{CMC and constant Gauss curvature foliations}

Let $M$ be a 3-dimensional globally hyperbolic AdS manifold. Recall that $M$ contains 
a unique minimal non-empty convex subset, called its convex core. It is known that:
\begin{itemize}
\item $M$ admits a unique foliation by closed, space-like, constant mean curvature surfaces, 
with the mean curvature going from $-\infty$ to $\infty$ (see \cite{BBZ}); in particular
$M$ contains a unique closed, space-like maximal surface;
\item the complement in $M$ of the convex core admits a unique foliation by constant Gauss
curvature surfaces, see \cite{BBZ2}.
\end{itemize}

Consider now the universal cover of $M$. It can considered as a convex domain 
$D\subset\AdS_3$, with boundary at infinity a weakly space-like curve $\Lambda$ in $\partial_\infty\AdS_3$ which
is the limit set of the action of $\pi_1(M)$ on $D$. In the natural identification of
$\partial_\infty\AdS_3$ with $\R P^1\times \R P^1$, the curve $\Lambda$ is the graph of the quasi-symmetric
homeomorphism from $\R P^1$ to $\R P^1$ which conjugates the left and right actions of $\pi_1(M)$
on $\R P^1=\partial_\infty\HH^2$, see \cite{mess}.
Let $C$ be the convex hull of $\Lambda$. It follows from the two facts recalled above that:
\begin{itemize}
\item $D$ has a foliation by space-like surfaces of constant curvature with boundary at infinity $\Lambda$
(and this foliation is unique among those invariant under the action of $\pi_1(M)$ on $D$).
\item $D\smallsetminus C$ has a foliation by space-like constant Gauss curvatures with boundary at infinity
$\Lambda$ (and it is unique among those foliations invariant under the action of $\pi_1(M)$ on $D$).
\end{itemize}

Now remove the invariance hypothesis, and let $\Lambda\subset \partial_\infty \AdS_3$ be the graph of any
quasi-symmetric homeomorphism from $\R P^1$ to $\R P^1$, let $C$ be the convex hull of $\Lambda$, and 
let $D$ be the domain of dependence defined by $\Lambda$ (\emph{i.e.}\ the set of points $x\in \AdS_3$ such that all
time-like geodesics through~$x$ intersect $C$). It is known (see \cite{maximal}) that there
is a unique maximal space-like surface in $D$ with bounded principal curvatures and boundary at infinity $\Lambda$.

\begin{question} \label{q:mean}
Is there a unique foliation of $D$ by space-like, constant mean curvature surfaces with bounded
principal curvature and boundary at infinity $\Lambda$?
\end{question}

\begin{question} \label{q:gauss}
Is there a unique foliation of $D\smallsetminus C$ by space-like, constant Gauss curvature surfaces
with bounded principal curvatures?
\end{question}

The main motivation for the result of \cite{maximal} on maximal surfaces is that it provides an existence
and uniqueness result for minimal Lagrangian extensions to the hyperbolic plane of quasi-symmetric homeomorphisms
of the circle. A positive answer to either Question \ref{q:mean} or Question \ref{q:gauss} might similarly
have a translation in terms of the existence of special extensions to the hyperbolic disk of quasi-symmetric homeomorphisms
of the circle.

\section{Massive particles and tachyons}\label{sec:tachyons}

As in Sections \ref{sec:convex} and~\ref{sec:volume}, we denote by $S$ a closed surface with genus $g\geq 2$ and by $\cT$ the Teichm\"uller space of $S$.
We set $M=S\times\R$ and denote by $\cGH$ the space of globally hyperbolic maximal compact AdS structures on $M$, considered up to isotopy.

\subsection{The moduli space of AdS spacetimes with interacting particles}

If a 3-dimensional manifold contains massive particles --- cone singularities along
time-like lines of angle less than $2\pi$ --- then those particles can interact: the 
singular set can be a graph, with vertices corresponding to interaction points. 
When no singularity is present, the moduli space of the corresponding geometric
structures is parameterized by the product of two copies of $\cT$, see \cite{mess},
and this extends to the situation where the particles have angles less than $\pi$
so that they cannot interact, see \cite{cone}.

Can these results be extended in some way to the situation where the angles are 
only less than $2\pi$ and interactions occur? A possible answer is proposed in 
\cite{colII}, with the parameterization data given as a sequence of pairs of
hyperbolic metrics with cone singularities, with simple rules on how to go
from one pair to the next. It is proved there that this data provides a 
local parameterization of the space of AdS structures with interacting
particles. It would be interesting to understand whether this parameterization
is global, or whether some additional condition on the parameterization data
must be imposed to obtain a global parameterization.

\subsection{Rigidity of AdS manifolds with tachyons}

We return to the case of GHMC AdS structures on $M = S \times \R$ without particles. 
Consider the following restriction of Question~\ref{q:lamination} to the case of \emph{rational bending measures}:

\begin{question} \label{q:multi-curves}
If $l_-$ and $l_+$ are \emph{weighted multi-curves} which fill $S$, is there a unique GHMC AdS structure
$h$ on $M$ such that $l_-$ and $l_+$ are the measured bending laminations of the boundary of the 
convex core?
\end{question}

In the hyperbolic quasi-Fuchsian case, uniqueness was shown by Bonahon--Otal \cite{bonahon-otal} 
by doubling the convex core to produce a compact cone manifold and then applying Hogdson--Kerckhoff \cite{HK} 
local rigidity theory and other tools. 
In the $\AdS$ case, doubling a convex core with bending along simple closed curves produces 
a compact AdS manifold with \emph{tachyon} singularities, as defined in \cite{colI}. 
A tachyon is a singularity along a space-like geodesic so that the holonomy around the singular locus is a Lorentz boost. The magnitude of the boost gives the (absolute value of the) \emph{mass} of the tachyon (note that tachyons coming from doubled convex cores always have negative mass). A tachyon singularity is exactly the Lorentzian analogue of a cone singularity in a Riemannian 3-manifold. Following the approach of Bonahon--Otal, one could ask for an AdS version of Hodgson--Kerckhoff local rigidity:

\begin{question}
Let $N$ be a closed AdS manifold with tachyons, so that the singular locus is a link. Is the AdS geometry of $N$ rigid relative to the tachyon mass?
\end{question}

Of course, the answer to this question is no in the most basic case: 
In contrast to hyperbolic structures, $\AdS$ structures without singularities on a closed 3-manifold are very flexible (see Section~\ref{subsec:compactAdS3}). 
Such manifolds are Seifert fibered, and one way to produce deformations of the $\AdS$ geometry is to deform the base of the fibration 
(Goldman \cite{Goldman.Non.Standard} gives a construction of deformations arising in a different way). 
However, examples in Danciger's thesis \cite{Danciger-11} suggest that the flexibility of AdS manifolds is often greatly constricted by the presence of singularities. In particular, when the complement of the singular locus is atoroidal, these manifolds seem to exhibit rigidity properties (rel boundary data) that mimic those of hyperbolic structures. So we refine our previous question:

\begin{question} \label{q:tachyons-atoroidal}
Let $N$ be a closed AdS manifold with tachyons, so that the singular locus is a link. Further, assume that the complement of the singular locus is atoroidal. Is the AdS geometry of~$N$ rigid relative to the tachyon mass?
\end{question}

Note that AdS manifolds gotten by doubling a GHMC convex core (with rational bending) indeed satisfy the atoroidal hypothesis. This follows from the filling property of the bending laminations.

\begin{remark}
A positive answer to Question~\ref{q:tachyons-atoroidal} is not required by the uniqueness in Question~\ref{q:multi-curves}. As far as we know, the following could happen: There may be non-trivial deformations of a doubled convex core with tachyon mass remaining constant, but so that the deformed manifold is no longer a double. 
\end{remark}

Many closed AdS manifolds with tachyon singularities satisfying the atoroidal hypothesis of Question~\ref{q:tachyons-atoroidal} are produced by the geometric transition technology introduced by Danciger \cite{Danciger-11}. Roughly, the construction produces AdS manifolds with tachyons as the continuation of a path of hyperbolic cone structures that collapse onto a two-dimensional surface. Such transitioning paths of geometric structures can also be constructed directly using ideal tetrahedra (again see \cite{Danciger-11}). The existence of a geometric transition between hyperbolic cone manifolds and AdS tachyon manifolds could be regarded as evidence toward an affirmative answer to Question~\ref{q:tachyons-atoroidal}. One possible approach to Question~\ref{q:tachyons-atoroidal} and other similar questions is to develop a more effective theory of geometric transitions that would allow one to translate properties of the hyperbolic manifolds on one side of the transition to the!
  AdS manifolds on the other side.

We conclude this section with one more question.
\begin{question} \label{q:small-mass}
Let $N$ be a closed 3-manifold with a non-singular AdS structure and let $\Sigma$ be a link in $N$. Do there exist nearby AdS structures with tachyon singularities along $\Sigma$? What topological conditions on $(N,\Sigma)$ are needed to have such deformations? Geometric conditions?
\end{question}

Depending on $N$ and $\Sigma$, tachyon singularities cannot always be introduced. For example, let $N$ be the unit tangent bundle of a hyperbolic surface $S$ and consider the standard AdS structure on $N$ induced by a hyperbolic structure on $S$. Now let $\Sigma$ be a circle whose projection down to $S$ is a separating curve $c$. It is an easy exercise to show that a tachyon singularity along $\Sigma$ cannot be introduced by a small deformation of the AdS geometry. 

\section{The complex AdS (or hyperbolic) space in dimension 3}\label{sec:complex}

Consider the complex symmetric bilinear form
$$b(z,z') = - z_1z'_1 - z_2z'_2 + z_3z'_3 + z_4z'_4$$
on~$\C^4$ and let $X_{\C}$ be the set of points $z\in\C^4$ such that $b(z,z)=-1$.
Then $X_{\C}$ identifies with the 3-dimensional complex sphere $\mathbb{S}^3_{\C}=\OO(4,\C)/\OO(3,\C)$; its isometry group $\OO(4,\C)$ is $\PSL_2(\C)\times\nolinebreak\PSL_2(\C)$ up to finite index and a $2$-fold covering.
The space $X_{\C}$ contains both $\AdS_3=X_{\C}\cap\R^4$ and $\HH^3=X_{\C}\cap(\R\times\R\times i\R\times\R^+)$, and is a complexification of both of them.

\subsection{Equivariant surfaces}

Let $S$ be a closed surface.

\begin{question}
Given two quasi-Fuchsian representations of $\pi_1(S)$ into $\PSL_2(\C)$, is there an equivariant space-like embedding
$\widetilde{S}\rightarrow X_{\C}$? In this case, may
we find such an embedding of constant curvature?
\end{question}

Here $\widetilde{S}$ denotes a universal covering of~$S$.
We say that an embedding $\widetilde{S}\rightarrow X_{\C}$ is \emph{space-like} if the restriction of the complex form~$b$ to the tangent spaces is positive definite.

\begin{question}\label{complexproper}
Given two quasi-Fuchsian representations $\rho_l,\rho_r$ of $\pi_1(S)$ into $\PSL_2(\C)$, is there a (maximal)
domain of $X_{\C}$ where the action of the group $(\rho_l,\rho_r)(\pi_1(S))\subset\PSL_2(\C)\times\PSL_2(\C)$ is properly discontinuous? What is the topology of the corresponding quotient?
\end{question} 

\subsection{Proper actions on $X_{\C}$}

In the setting of Question~\ref{complexproper}, the action of $(\rho_l,\rho_r)(\pi_1(S))$ can never be properly discontinuous on \emph{the whole of~$X_{\C}$}.
However, if we allow $S$ to have boundary, then there exist pairs $(\rho_l,\rho_r)$ of convex cocompact representations of $\pi_1(S)$ into $\PSL_2(\C)$ such that the action on~$X_{\C}$ is proper.
This follows from \cite{gueritaud-kassel}, where the theory of Sections \ref{subsec:compactAdS3} and~\ref{subsec:coreproper} is generalized in the following way: given any hyperbolic $3$-manifold~$M$ and any pair of convex cocompact representations $\rho_l,\rho_r : \pi_1(M)\rightarrow\PSL_2(\C)$, the action of the group $(\rho_l,\rho_r)(\pi_1(M))$ on~$X_{\C}$ is properly discontinuous if and only if 
$C_{length}(\rho_l,\rho_r)<1$, or equivalently if and only if $C_{Lip}(\rho_l,\rho_r)<\nolinebreak 1$, where $C_{length}$ and $C_{Lip}$ are defined similarly to Sections \ref{subsec:compactAdS3} and~\ref{subsec:coreproper}.
One can ask an analogue of Questions \ref{readproperness} and~\ref{complexproper} in this context.

\begin{question}
Let $M$ be a hyperbolic $3$-manifold and let $(\rho_l,\rho_r)$ be a pair of convex cocompact representations of $\pi_1(M)$ into $\PSL_2(\C)$.
Is there a (maximal) domain of~$X_{\C}$ where the action of the group $(\rho_l,\rho_r)(\pi_1(M))\subset\PSL_2(\C)\times\PSL_2(\C)$ is properly discontinuous? What is the topology of the corresponding quotient?
How can one read geometrically the properness or non-properness of the action of $(\rho_l,\rho_r)(\pi_1(M))\subset\PSL_2(\C)\times\PSL_2(\C)$ on~$X_{\C}$?
\end{question}

We note that $X_{\C}$ is the model space for all \emph{holomorphic Riemannian} 
complex $3$-manifolds of constant non-zero curvature.
The importance of these manifolds is highlighted by the uniformization theorem of Dumitrescu--Zeghib 
\cite{dumitrescu-zeghib}, which states that any $3$-dimensional compact complex manifold~$M$ admitting a holomorphic 
Riemannian metric admits one of constant sectional curvature, up to a finite covering.
For such a constant-curvature metric, $M$ should conjecturally always be complete, \emph{i.e.}\ realized (when the curvature is non-zero) as a quotient of~$X_{\C}$ by a discrete subgroup of $\PSL_2(\C)\times\nolinebreak\PSL_2(\C)$ 
acting properly discontinuously \cite{dumitrescu-zeghib}.


\newcommand{\etalchar}[1]{$^{#1}$}
\def\cprime{$'$}

\end{document}